\documentclass[a4paper]{amsart}

\usepackage{amssymb}
\usepackage{amscd}
\usepackage{amsthm}
\usepackage{amsmath}
\usepackage{latexsym}
\usepackage[all]{xy}
\usepackage{amsfonts}

\usepackage{booktabs}
\usepackage{multirow}
\usepackage{colortbl}

\usepackage{hyperref}
\usepackage{extarrows}
\usepackage{graphicx}
\usepackage{tikz}
\usepackage[table]{}

\theoremstyle{plain}
\newtheorem{theorem}{Theorem}
\newtheorem{corollary}[theorem]{Corollary}
\newtheorem{lemma}[theorem]{Lemma}

\theoremstyle{definition}
\newtheorem{definition}[theorem]{Definition}
\newtheorem{example}[theorem]{Example}

\theoremstyle{remark}

\newtheorem{remark}[theorem]{Remark}

\newtheorem{acknowledgment}[theorem]{Acknowledgment}
\numberwithin{theorem}{section}

\newcommand{\card}[1]{\mbox{\rm{card}}(#1)}

\newcommand{\Soc}[1]{\mbox{\rm{Soc}}(#1)}

\newcommand{\im}[1]{\mbox{\rm{Im}}(#1)}

\newcommand{\Hom}[3]{\mbox{\rm{Hom}}_{#1}(#2,#3)}
\newcommand{\Ext}[4]{\mbox{\rm{Ext}}^{#1}_{#2}(#3,#4)}

\newcommand{\rmod}[1]{\mbox{\rm{Mod}--}{#1}}

\newcommand{\Ker}[1]{\mbox{\rm{Ker}}(#1)}

\begin{document}

\title[Weak diamond, weak projectivity, and transfinite extensions]{Weak diamond, weak projectivity, and transfinite extensions of simple artinian rings}
\author{\textsc{Jan Trlifaj}}
\address{Charles University, Faculty of Mathematics and Physics, Department of Algebra \\
Sokolovsk\'{a} 83, 186 75 Prague 8, Czech Republic}
\email{trlifaj@karlin.mff.cuni.cz}

\date{\today}
\thanks{Research supported by GA\v CR 20-13778S}
\begin{abstract} We apply set-theoretic methods to study projective modules and their generalizations over transfinite extensions of simple artinian rings $R$. We prove that if $R$ is small, then the Weak Diamond implies that projectivity of an arbitrary module can be tested at the layer epimorphisms of $R$.        
\end{abstract}

\maketitle

The classic Baer's Criterion, saying that a module $M$ is injective, iff it is $R$-injective, is a basic tool of the structure theory of injective modules over an arbitrary ring $R$. However, unless $R$ is a perfect ring, there are no criteria available for the dual case, that is, for testing projectivity using a set of epimorphisms, \cite{ST}. 

For many non-perfect rings, it can be proved that there exist small (e.g., countably generated) non-projective $R$-projective modules (e.g., when $R$ is commutative noetherian of Krull dimension $\geq 1$). However, for each cardinal $\kappa$, there exists a non-right perfect ring $R_\kappa$ such that all $\leq \kappa$-generated $R_\kappa$-projective modules are projective, \cite{T2}. This is the best one can achieve in ZFC, because it is consistent with ZFC + GCH that if $R$ is not right perfect, then there always exist (large) $R$-projective modules that are not projective, cf.\ \cite{AIPY}.

Consistency (and hence independence) of the coincidence of $R$-projectivity and projectivity for certain commutative non-noetherian rings was proved in \cite{T1}. This answered in the positive a question from \cite[2.8]{AIPY}, and clarified the set-theoretic status of an old problem by Carl Faith \cite[p.175]{F}. The consistency result was extended in \cite{T2} to further classes of rings that are finite Loewy length extensions of simple artinian rings. The set theoretic tool used in \cite{T1} and \cite{T2} was Jensen's Diamond. 

\medskip
The goal of the present paper is twofold: to enhance the algebraic tools to cover infinite Loewy length extensions of simple artinian rings, and to weaken the set-theoretic assumptions used in the proofs. In Theorem \ref{weak} below, we show that the Weak Diamond Principle $\Phi$ and CH are sufficient to prove coincidence of the classes of all weak $R$-projective, $R$-projective, and projective modules, in the case when $R$ has cardinality at most $\aleph_1$, its Loewy length is countable, and each proper layer of $R$ is countably generated. That is, in this case, $\Phi$ and CH imply that the projectivity of a module $M$ is equivalent to the factorization of all morphisms from $M$ with finitely generated images through the layer epimorphisms of $R$. The latter are just the canonical projections $\pi_\alpha : S_{\alpha + 1} \to S_{\alpha + 1}/S_\alpha$ ($\alpha < \sigma$), where $(S_\alpha \mid \alpha \leq \sigma + 1)$ is the socle sequence of $R$.  

\medskip  
For basic notions and facts needed from ring and module theory, we refer to \cite{AF} and \cite{G}; our references for set-theoretic homological algebra are \cite{EM} and \cite{GT}.          

\section{Weak projectivity}

We start by recalling the classic notion of relative projectivity from \cite[\S16]{AF}: 

\begin{definition}\label{Nproj} Let $R$ be a ring and $M$, $N$ be modules. Then $M$ is \emph{$N$-projective} provided that for each submodule $P$ of $N$,  each homomorphism $f : M \to N/P$ has a factorization through the canonical projection $N \to N/P$. That is, the functor 
$\Hom RM{-}$ is exact on all short exact sequences whose middle term is $N$.        
\end{definition}

The following Lemma is well-known (see \cite[16.12 and 16.14]{AF}):

\begin{lemma}\label{af} The class of all modules $N$ such that $M$ is $N$-projective is closed under submodules, homomorphic images, and finite direct sums. 

In particular, each finitely generated $R$-projective module is projective. 
\end{lemma}

\medskip
Let $R$ be a right semiartinian ring with the right socle sequence $(S_\alpha \mid \alpha \leq \sigma + 1)$. We will call $\sigma + 1$ the \emph{Loewy length} of $R$. To avoid trivialities, we will tacitly assume that $\sigma > 0$, that is, that $R$ is not completely reducible. 

For each $\alpha \leq \sigma$, we will call the completely reducible module $L_\alpha = S_{\alpha +1}/S_\alpha$ the \emph{$\alpha$th layer} of $R$. The layers $L_\alpha$ for $0 < \alpha < \sigma$ are called \emph{proper}. For each $\alpha \leq \sigma$, the canonical epimorphism $\pi_\alpha : S_{\alpha +1} \to L_\alpha$ is the \emph{$\alpha$th layer epimorphism} of $R$.   

In this setting, the following definition of weak $R$-projectivity was introduced in \cite[3.6]{T2} 

\begin{definition}\label{weakproj} A module $M$ is called \emph{weakly R-projective} provided that for each $0 < \alpha \leq \sigma$, each $f\in \Hom RM{L_\alpha}$ with a finitely generated image has a factorization through the $\alpha$th layer epimorphism $\pi_\alpha$.
\end{definition}

Note that by Lemma \ref{af}, each $R$-projective module is weakly $R$-projective. The following easy observation \cite[2.3]{T2} shows that weak $R$-projectivity can equivalently be stated in a stronger form:

\begin{lemma}\label{stronger} Let $R$ be a right semiartinian ring with the right socle sequence $(S_\alpha \mid \alpha \leq \sigma + 1)$. Let $M$ be a module. Then the following are equivalent:
\begin{enumerate}
\item $M$ is weakly $R$-projective, 
\item For each $\alpha \leq \sigma$, each submodule $K$ such that $S_\alpha \subseteq K \subseteq S_{\alpha + 1}$ and each 
$f\in \Hom RM{S_{\alpha +1}/K}$ with a finitely generated image, there exists $g \in \Hom RM{S_{\alpha +1}}$ such that $f = \pi_K g$, where $\pi_K : S_{\alpha +1} \to S_{\alpha +1}/K$ is the canonical projection.  
\end{enumerate}
\end{lemma}
  
In the proof of the next lemma, we will verify $N$-projectivity of a finitely generated module $M$ by a recursive procedure using a filtration of $N$. 

Here, we call a chain $( N_\alpha \mid \alpha \leq \tau )$ of submodules of a module $N$ a \emph{filtration} of $N$, provided that $\tau$ is an ordinal, $N_0 = 0$, $N_\alpha \subseteq N_{\alpha +1}$ for each $\alpha < \tau$, $N_\alpha = \bigcup_{\beta < \alpha} N_\beta$ for each limit ordinal $\alpha \leq \tau$, and $N_\tau = N$.     

\begin{lemma}\label{finNproj} Let $R$ be a ring and $N$ be a module with a filtration $\mathcal N = ( N_\alpha \mid \alpha \leq \tau)$.  

Let $M$ be a finitely generated module such that for each $\alpha < \tau$, each submodule $K$ such that $N_\alpha \subseteq K \subseteq N_{\alpha + 1}$ and each $f\in \Hom RM{N_{\alpha +1}/K}$, there exists $g \in \Hom RM{N_{\alpha +1}}$ such that $f = \pi_K g$, where $\pi_K : N_{\alpha +1} \to N_{\alpha +1}/K$ is the canonical projection. 

Then $M$ is $N$-projective.
\end{lemma}
 \begin{proof} We claim that $M$ is $N_\alpha$-projective for each $\alpha \leq \tau$. If so, then for $\alpha = \tau$, we get that $M$ is $N$-projective. 

The claim will be proved by induction on $\alpha$. There is nothing to prove for $\alpha = 0$.

Assume that the claim is true for some $\alpha < \tau$. Let $K$ be a submodule of $N_{\alpha + 1}$. Let $\pi _K : N_{\alpha + 1} \to N_{\alpha + 1}/K$, $\rho : N_{\alpha + 1} \to N_{\alpha + 1}/(N_\alpha + K)$, $\eta : N_{\alpha + 1}/K \to N_{\alpha + 1}/(N_\alpha + K)$, and $\theta : N_\alpha \to N_\alpha /(N_\alpha \cap K)$ denote the canonical projections. Also, let $\iota _\alpha : N_\alpha / (N_\alpha \cap K) \to (N_\alpha + K)/K$ be the canonical isomorphism (given by the assignment $\iota _\alpha (x + N_\alpha \cap K) = x + K$).

These homomorphisms fit in the following commutative diagram with exact rows (where, except for the zero maps, all the unnamed single arrows are inclusions):

$$\begin{CD}
0 @>>> N_\alpha + K @>>> N_{\alpha + 1} @>{\rho}>> N_{\alpha + 1}/(N_\alpha + K) @>>> 0 \\
@. @AAA @| @A{\eta}AA @. \\
0 @>>> K @>>> N_{\alpha + 1} @>{\pi_K}>> {N_{\alpha + 1}/K} @>>> 0 \\
@. @| @AAA @AAA @. \\
0 @>>> K @>>> {N_\alpha + K} @>{\pi_K \restriction}>> {(N_\alpha + K)/K} @>>> 0 \\
@. @AAA  @AAA @A{\iota _\alpha}AA @. \\
0 @>>> {N_\alpha \cap K} @>>> {N_\alpha} @>{\theta}>> {N_\alpha/(N_\alpha \cap K)} @>>> 0
\end{CD}$$

Let $f \in \Hom RM{N_{\alpha + 1}/K}$. We have to show that $f$ factorizes through $\pi_K$.

By our assumption on $M$, there exists $g \in \Hom RM{N_{\alpha +1}}$ such that $\rho g = \eta f$. Since $\rho = \eta \pi _K$, $\eta (f - \pi _K g) = 0$. It follows that $\delta = f - \pi _K g$ maps $M$ into $\ker (\eta) = (N_\alpha + K)/K$, whence  $\iota_{\alpha}^{-1} \delta \in \Hom RM{N_{\alpha}/(N_\alpha \cap K)}$. By the inductive premise, there exists $\epsilon : M \to N_\alpha$ such that $\iota _\alpha \theta \epsilon = \delta$. As $\iota _\alpha \theta = \pi _K \restriction (N_\alpha + K)$, we conclude that $f = \pi _K g + \delta = \pi _K g + \pi _K \epsilon = \pi _K (g + \epsilon)$, which is the desired factorization of $f$ through $\pi _K$.

Let $\alpha \leq \tau$ be a limit ordinal.  Let $K$ be a submodule of $N_{\alpha} = \bigcup_{\beta < \alpha} N_\beta$ and $f \in \Hom RM{N_{\alpha}/K}$. Since $M$ is finitely generated, there exists $\beta < \alpha$ such that $f : M \to (N_{\beta} + K)/K$. Let $\iota _\beta : N_{\beta}/(N_\beta \cap K) \to (N_\beta + K)/K$ be the canonical isomorphism, and $\theta : N_\beta \to N_\beta / (N_\beta \cap K)$ the canonical projection. Then $\pi _K \restriction (N_\beta + K) = \iota _\beta \theta$, and we have the following commutative diagram with exact rows (again, except for the zero maps, all the unnamed single arrows are inclusions):

$$\begin{CD}
0 @>>> K @>>> N_{\alpha} @>{\pi_K}>> {N_{\alpha}/K} @>>> 0 \\
@. @| @AAA @AAA @. \\
0 @>>> K @>>> {N_\beta + K} @>{\pi_K \restriction}>> {(N_\beta + K)/K} @>>> 0 \\
@. @AAA  @AAA @A{\iota _\beta}AA @. \\
0 @>>> {N_\beta \cap K} @>>> {N_\beta} @>{\theta}>> {N_\beta/(N_\beta \cap K)} @>>> 0
\end{CD}$$

By the inductive premise for $\beta$, there exists $g \in \Hom RM{N_\beta}$ such that $\theta g = \iota _{\beta}^{-1} f$. Then $f = \iota _{\beta} \theta g = \pi _K g$, which yields the desired factorization of $f$ through $\pi _K$. 
\end{proof}

The results above imply that weak $R$-projectivity is sufficient to guarantee projectivity for each finitely generated module:  
 
\begin{corollary}\label{finproj} Let $R$ be any right semiartinian ring. Let $M$ be a finitely generated weakly $R$-projective module. Then $M$ is projective.
\end{corollary}
\begin{proof} Since weak projectivity can be expressed in the stronger form of Lemma \ref{stronger}, we can apply Lemma \ref{finNproj} for $N = R$ and $\mathcal N$ = the right socle sequence of $R$. Thus, $M$ is $R$-projective, and by Lemma \ref{af}, $M$ is projective.
\end{proof}

\medskip
The question of the coincidence of weak $R$-projectivity, $R$-projectivity, and projectivity for infinitely generated modules over semiartinian rings, that is, whether projectivity can be tested at the layer epimorphisms of $R$, is much more delicate. We will see that for some semiartinian rings, it is actually independent of ZFC + GCH. 

On the one hand, there is the following set-theoretic barrier for all non-right perfect rings $R$, given by Shelah's Uniformization Principle (SUP). This principle is consistent with ZFC + GCH (see \cite[XIII.1.5]{EM} or \cite[\S2]{ES} for more details):  

\begin{lemma}\label{sup} Assume SUP. Let $R$ be a non-right perfect ring. Let $\kappa$ be a singular cardinal of cofinality $\omega$ such that $\card R < \kappa$. Then there exists a $\kappa^+$ generated module $M$ of projective dimension equal to $1$, such that $\Ext 1RMI = 0$ for each right ideal $I$ of $R$. In particular, $M$ is $R$-projective, but not projective.
\end{lemma}  
\begin{proof} See \cite[2.4 and 2.5]{AIPY}.
\end{proof}    

On the other hand, there are different extensions of ZFC + GCH where weak $R$-projectivity, $R$-projectivity, and projectivity coincide for suitable semiartinian rings. These rings will be studied in the next section.  

\section{Transfinite extensions of simple artinian rings}

Recall that a ring $R$ is \emph{von Neumann regular}, if for each $r \in R$ there exists $s \in R$ such that $rsr = r$. We refer to \cite{G} for properties of von Neumann regular rings. A~ring $R$ is said to have (right) \emph{primitive factors artinian}, or (right) \emph{pfa} for short, if $R/P$ is (right) artinian for each (right) primitive ideal of $R$.  

If $R$ is von Neumann regular, then minimal right (left) ideals of $R$ correspond to primitive idempotents of $R$, so the right and left socle sequences of $R$ coincide. In particular, $R$ is right semiartinian, iff $R$ is left semiartinian. Similarly, by \cite[6.2]{G}, a von Neumann regular ring has right pfa, iff it has left pfa. The latter condition can equivalently be stated as a property of the module category $\rmod R$: each homogenous completely reducible module is injective, see \cite[6.28]{G}.             

\medskip 
Let $R$ be a semiartinian von Neumann regular ring with pfa, and $(S_\alpha \mid \alpha \leq \sigma + 1)$ be its socle sequence.  (Notice that $\sigma > 0$ implies that $R$ is not right perfect, because the Jacobson radical of any von Neumann regular ring is $0$; thus, Lemma \ref{sup} applies here.) 

By the following theorem from \cite{RTZ}, semiartinian von Neumann regular rings with pfa can be viewed as transfinite extensions of full matrix rings over skew-fields (i.e., of simple artinian rings):

\begin{theorem}\label{RTZth} Let $R$ be a right semiartinian ring and $(S_\alpha \mid \alpha \leq \sigma + 1)$ be the right socle sequence of $R$. The following conditions are equivalent: 
\begin{enumerate}
\item $R$ is von Neumann regular with pfa.
\item for each $\alpha \leq \sigma$ there are a cardinal $\lambda_\alpha$, positive integers $n_{\alpha\beta}$ ($\beta < \lambda_\alpha$) and skew-fields $K_{\alpha\beta}$ ($\beta < \lambda_\alpha$), such that $L_\alpha \overset{\varphi_\alpha}\cong \bigoplus_{\beta < \lambda_\alpha} M_{n_{\alpha\beta}}(K_{\alpha\beta})$, as rings without unit. Moreover, $\lambda_\alpha$ is infinite iff $\alpha < \sigma$. 

The pre-image of $M_{n_{\alpha\beta}}(K_{\alpha\beta})$ in the isomorphism $\varphi_\alpha$ coincides with the $\beta$th homogenous component of $\Soc{R/S_\alpha}$, and it is a finitely generated as right $R/S_\alpha$-module for all $\beta < \lambda_\alpha$. 
\end{enumerate} 
\end{theorem}

The structure of the rings characterized by Theorem \ref{RTZth} can be depicted as follows:

\eject

\medskip

\begin{table}[!ht]
\begin{tabular}{ll|llll|l}
\cline{3-6}
\cellcolor[HTML]{34CDF9}{$\mathbf{L_\sigma}$}  &  & \cellcolor[HTML]{34CDF9}{$\mathbf{M_{n_{\sigma 0}}(K_{\sigma 0}) \, \oplus}$} & \cellcolor[HTML]{34CDF9}{...} & \cellcolor[HTML]{34CDF9}{$\mathbf{\oplus \, M_{n_{\sigma , \lambda_\sigma -1}}(K_{\sigma, \lambda_\sigma -1})}$} &                                                           &                          \\ \cline{3-6}
\cellcolor[HTML]{96FFFB}{...} &  & \cellcolor[HTML]{96FFFB}{...} & \cellcolor[HTML]{96FFFB}{...} & \cellcolor[HTML]{96FFFB}{...} & \multicolumn{1}{l|}{\cellcolor[HTML]{96FFFB}\textbf{}} &                           \\ \cline{3-6}
\cellcolor[HTML]{38FFF8}{$\mathbf{L_\alpha}$}  &  & \cellcolor[HTML]{38FFF8}{$\mathbf{M_{n_{\alpha 0}}(K_{\alpha 0}) \, \oplus}$} & \cellcolor[HTML]{38FFF8}{...}  & \cellcolor[HTML]{38FFF8}{$\mathbf{\oplus \, M_{n_{\alpha\beta}}(K_{\alpha \beta}) \, \oplus}$ ...} & \multicolumn{1}{l|}{\cellcolor[HTML]{38FFF8}{}} & $\beta < \lambda_\alpha$          \\ \cline{3-6}
\cellcolor[HTML]{68CBD0}{...} &  & \cellcolor[HTML]{68CBD0}{...} & \cellcolor[HTML]{68CBD0}{...} & \cellcolor[HTML]{68CBD0}{...} & \multicolumn{1}{l|}{\cellcolor[HTML]{68CBD0}\textbf{}} &                           \\ \cline{3-6}
\cellcolor[HTML]{00D2CB}{$\mathbf{L_1}$}  &  & \cellcolor[HTML]{00D2CB}{$\mathbf{M_{n_{10}}(K_{10}) \, \oplus}$} & \cellcolor[HTML]{00D2CB}{...} & \cellcolor[HTML]{00D2CB}{$\mathbf{\oplus \, M_{n_{1\beta}}(K_{1 \beta}) \, \oplus}$ ...} & \multicolumn{1}{l|}{\cellcolor[HTML]{00D2CB}\textbf{}} & $\beta < \lambda_1$          \\ \cline{3-6}
\cellcolor[HTML]{329A9D}{$\mathbf{L_0 = Soc(R)}$}  &  & \cellcolor[HTML]{329A9D}{$\mathbf{M_{n_{00}}(K_{00}) \, \oplus}$}  & \cellcolor[HTML]{329A9D}{...} & \cellcolor[HTML]{329A9D}{$\mathbf{\oplus \, M_{n_{0 \beta}}(K_{0 \beta}) \, \oplus}$ ...} & \multicolumn{1}{l|}{\cellcolor[HTML]{329A9D}\textbf{}} & $\beta < \lambda_0$ \\ \cline{3-6}
\end{tabular}
\end{table}

\medskip

Here, $\sigma + 1$ is the Loewy length of $R$. The rows in this picture represent the layers of $R$; $\lambda_\alpha$ is the number of homogenous components in the $\alpha$th layer $L_\alpha = S_{\alpha + 1}/S_\alpha$ for each $\alpha \leq \sigma$ (this number is infinite except for $\alpha = \sigma$). 

$n_{\alpha\beta}$ is the (finite) rank of the $\beta$th homogenous component of $L_\alpha$, and $K_{\alpha\beta}$ is the endomorphism ring (skew-field) of each simple module in the $\beta$th homogenous component of $L_\alpha$ (for all $\alpha \leq \sigma$, $\beta < \lambda_\alpha$). 

In particular, $\sigma$, $\lambda_\alpha$ ($\alpha \leq \sigma$), $n_{\alpha\beta}$, and $K_{\alpha\beta}$ ($\alpha \leq \sigma$, $\beta < \lambda_\alpha$) are invariants of the ring $R$. 

\begin{lemma}\label{count} Let $R$ be a semiartinian von Neumann regular ring with pfa. Then the following hold:
\begin{enumerate}
\item The class of all weakly $R$-projective modules is closed under submodules. 
\item All countably generated weakly $R$-projective modules are projective. 
\end{enumerate}
\end{lemma}
\begin{proof}
(1) This follows from the fact \cite[6.28]{G} that all homogenous semisimple modules are injective, whence so are all finitely generated semisimple modules, and thus all finitely generated submodules of each of the layers $L_\alpha$ ($\alpha < \sigma$).

(2) Let $C$ be a countably generated weakly $R$-projective module. So $C = \bigcup_{n < \omega} F_n$ where $( F_n \mid n < \omega )$ is a chain of finitely generated submodules of $C$. By part (1), each $F_n$ is weakly $R$-projective, hence projective by Corollary \ref{finproj}. Since $R$ is von Neumann regular, $F_n$ is a direct summand in $F_{n+1}$ for each $n < \omega$, whence $C$ is projective.   
\end{proof}

From Lemma \ref{count}(1), we see that if the classes of all weak $R$-projective and projective modules coincide, then $R$ is a right hereditary ring. What that means in out setting is partially clarified in the next lemma:

\begin{lemma}\label{hereditary} Let $R$ be a semiartinian von Neumann regular ring with pfa. Consider the following two conditions:
\begin{enumerate}
\item $R$ is (left and right) hereditary.
\item $\sigma$ is a countable ordinal, and all proper layers of $R$ are countably generated (i.e., $\lambda_\alpha$ is countable for each $0 < \alpha < \sigma$).
\end{enumerate}
Then (2) implies (1). Moreover, if $\lambda_0$ is countable, then (1) and (2) are equivalent.
\end{lemma}
\begin{proof} That (2) implies (1) was proved in \cite[3.10]{T2}. 

Assume $\lambda_0$ is countable and $R$ is right hereditary. Since all the homogenous components of $\Soc{R}$ are injective, and hence finitely generated, $\Soc{R}$ is countably generated. 

Assume (2) fails. If $\sigma$ is uncountable, then $\sigma = \tau + n$ for an uncountable limit ordinal $\tau$ and some $n < \omega$, whence $I = S_\tau$ is not countably generated. If $\lambda_\alpha$ is uncountable for some ordinal $0 < \alpha < \sigma$, then the $\alpha$th layer $L_\alpha$ is not countably generated, and the same is true of $I = S_{\alpha + 1}$. 

Since projective modules over von Neumann regular rings are isomorphic to direct sums of cyclic modules generated by idempotents of $R$, in either case there is a cardinal $\kappa > \aleph_0$ such that $I \cong \bigoplus_{\gamma < \kappa} e_{\gamma} R$ where $e_\gamma ^ 2 = e_\gamma \in R$ for each $\gamma < \kappa$. Since $R$ is semiartinian, $\Soc{I} \cong \bigoplus_{\gamma < \kappa} \Soc{e_{\gamma} R}$ is a direct sum of uncountably many non-zero completely reducible modules, in contradiction with $\Soc{I}$ being a direct summand in the countably generated module $\Soc{R}$.         
\end{proof}

\medskip 
Notice that Lemma \ref{hereditary} implies that all semiartinian von Neumann regular rings with pfa of Loewy length $2$ are hereditary. The simplest such example is the $K$-subalgebra $R_1$ of all eventually constant sequences in $K^\omega$ for a field $K$, studied in \cite{T1}. The corresponding picture for $R_1$ is a follows: 

\medskip

\begin{table}[!ht]
\begin{tabular}{
>{\columncolor[HTML]{34CDF9}}l l
>{\columncolor[HTML]{FFFFFF}}l 
>{\columncolor[HTML]{FFFFFF}}l 
>{\columncolor[HTML]{34CDF9}}c 
>{\columncolor[HTML]{FFFFFF}}l 
>{\columncolor[HTML]{FFFFFF}}l }
{$\mathbf{L_1 = K}$}                               &  & \textbf{}                                & \textbf{}                            & \textbf{K}                                     &                                      &                                      \\
\cellcolor[HTML]{96FFFB}{$\mathbf{L_0 = K^{(\omega)}}$} &  & \cellcolor[HTML]{96FFFB}{$\mathbf{K \, \oplus}$} & \cellcolor[HTML]{96FFFB}\textbf{...} & \cellcolor[HTML]{96FFFB}{$\mathbf{\oplus \, K \, \oplus}$} & \cellcolor[HTML]{96FFFB}\textbf{...} & \cellcolor[HTML]{96FFFB}\textbf{...}
\end{tabular}
\end{table}

\medskip

However, semiartinian von Neumann regular rings with pfa of Loewy length $3$ need not be hereditary: for an example, take a set $S$ of cardinality $\omega_1$ consisting of almost disjoint infinite subsets of $\omega$, and consider the $K$-subalgebra $R^\prime$ of $K^\omega$ generated (as a $K$-linear space) by $B \cup C \cup \{ 1 \}$ where $1$ is the unit of $K^\omega$, $B$ is the canonical basis of $K^{(\omega)}$ and $C$ is the set of characteristic functions of all the sets in $S$. By Lemma \ref{hereditary}, $R^\prime$ is not hereditary. The corresponding picture for $R^\prime$ is           
       
\medskip

\begin{table}[!ht]
\begin{tabular}{|l|lllllll}
\cline{1-1}
\cellcolor[HTML]{34CDF9}{$\mathbf{L_2 = K}$}              &  & \cellcolor[HTML]{FFFFFF}\textbf{}        & \cellcolor[HTML]{FFFFFF}\textbf{}    & \cellcolor[HTML]{34CDF9}\textbf{\,\, K}             & \cellcolor[HTML]{FFFFFF}             &                                      &                                      \\ 
\cline{1-1}
\cellcolor[HTML]{96FFFB}{$\mathbf{L_1 = K^{(\omega_{1})}}$} &  & \cellcolor[HTML]{96FFFB}{$\mathbf{K \, \oplus}$}  & \cellcolor[HTML]{96FFFB}\textbf{...} & \cellcolor[HTML]{96FFFB}{$\mathbf{\oplus \, K \, \oplus}$} & \cellcolor[HTML]{96FFFB}\textbf{...} & \cellcolor[HTML]{96FFFB}\textbf{...} & \cellcolor[HTML]{96FFFB}\textbf{...} \\ 
\cline{1-1}
\cellcolor[HTML]{38FFF8}{$\mathbf{L_0 = K^{(\omega)}}$}   &  & \cellcolor[HTML]{38FFF8}{$\mathbf{K \, \oplus}$}  & \cellcolor[HTML]{38FFF8}\textbf{...} & \cellcolor[HTML]{38FFF8}{$\mathbf{\oplus \, K \, \oplus}$} & \cellcolor[HTML]{38FFF8}\textbf{...} &                                      &                                      \\ 
\cline{1-1}
\end{tabular}
\end{table}

\medskip
The following recursive construction shows that hereditary semiartinian von Neumann regular rings with pfa of Loewy length $\alpha + 1$ do exist for each countable ordinal $\alpha > 0$. The induction step makes use of a construction of semiartinian von Neumann regular rings from \cite[2.4]{EGT}:

\begin{example}\label{exegt} Let $K$ be a field. By induction on $0 < \alpha < \omega_1$, we will construct semiartinian von Neumann regular $K$-algebras with pfa, $R_\alpha$, of Loewy length $\alpha +1$ such that $R_\alpha$ has countably generated layers, together with $K$-algebra embeddings $\nu_\alpha : R_\alpha \to K^\omega$, and for each $0 < \beta < \alpha$, non-unital $K$-algebra monomorphisms $f_{\alpha \beta}: R_\beta \to R_\alpha$ and $g_{\alpha \beta} : K^\omega \to K^\omega$ such that the squares 

\medskip
$$\begin{CD}
R_{\beta} @>{f_{\alpha \beta}}>> {R_{\alpha}}\\
@V{\nu_\beta}VV @V{\nu_{\alpha}}VV \\
{K^\omega} @>{g_{\alpha \beta}}>> {K^\omega}\\
\end{CD}$$

\medskip
are commutative, $g_{\alpha \beta}$ splits, and the complement of $\im{g_{\alpha \beta}}$ in $K^\omega$ is isomorphic to $K^\omega$. 
    
For $\alpha = 1$, we let $R_1$ be the $K$-algebra of all eventually constant sequences of elements of $K$ mentioned above. In particular, $R_1$ is a $K$-subalgebra of $K^\omega$ of Loewy length $2$, and $\nu_1$ is defined as the inclusion of $R_1$ into $K^\omega$.  

The induction step is modeled on \cite[2.4]{EGT}: Assume the construction is done up to some $0 < \alpha < \aleph_1$. Let $I = \bigoplus_{i < \omega} R_\alpha$. Then $\iota_\alpha = \bigoplus_{i<\omega} \nu_\alpha$ embeds $I$ into $D = \bigoplus_{i < \omega} K^\omega \subseteq K^{\omega^\omega}$. As in \cite{EGT}, $\omega^\omega$ denotes the ordinal exponentiation, so $\omega^\omega = \sup_{i < \omega} \omega^i$ and $\omega^{i+1} = \omega^i \times \omega$ for each $0 < i < \omega$. Note that $\iota_\alpha(I)$ is an ideal of the $K$-algebra $K^{\omega^\omega}$, and $S_{\alpha + 1} = \iota_\alpha(I) \oplus 1.K$ is a $K$-subalgebra of $K^{\omega^\omega}$. 

Since the ordinal $\omega^\omega$ is countable, there is a $K$-algebra isomorphism $\psi : K^{\omega^\omega} \to K^\omega$. Let $R_{\alpha + 1} = \psi (S_{\alpha + 1})$, and let $\nu_{\alpha + 1}$ denote the inclusion of $R_{\alpha + 1}$ into $K^\omega$. 

Let $\mu_{\alpha}$ be the embedding of the first copy of $K^\omega$ in $D$ composed with the inclusion $D \subseteq K^{\omega^\omega}$. Notice that $\mu_{\alpha}$ is a split non-unital embedding of $K$-algebras, and so is $\psi \mu_{\alpha}$. In fact, $\psi \mu_{\alpha}(K^\omega)$ is a direct summand in $K^{\omega}$ with a complement isomorphic to $K^{\omega}$. Moreover, $\mu_{\alpha} \restriction R_\alpha$ is a non-unital $K$-algebra embedding of $R_\alpha$ into $S_{\alpha + 1}$. 

We have the following commutative diagram, where the vertical maps are $K$-algebra embeddings, and the horizontal ones are non-unital $K$-algebra embeddings:

\medskip
$$\begin{CD}
R_{\alpha} @>{{\iota_{\alpha}}\restriction R_\alpha}>> {S_{\alpha + 1}} @>{\psi \restriction S_{\alpha +1}}>> {R_{\alpha + 1}}\\
@V{\nu_\alpha}VV @V{\subseteq}VV @V{\nu_{\alpha +1}}VV \\
{K^\omega} @>{\mu_{\alpha}}>> {K^{\omega^\omega}} @>{\psi}>> {K^\omega}\\
\end{CD}$$

\medskip
We let $f_{\alpha+1, \alpha} =  \psi \iota_\alpha \restriction R_{\alpha}$, $g_{\alpha+1, \alpha} =  \psi \mu_\alpha$, and for each $0 < \beta < \alpha$, $f_{\alpha + 1, \beta} = f_{\alpha+1, \alpha} f_{\alpha \beta}$ and $g_{\alpha + 1, \beta} = g_{\alpha+1, \alpha} g_{\alpha \beta}$

By \cite[2.4]{EGT}, the Loewy length of $R_{\alpha + 1}$ is $\alpha + 2$. Moreover, for each $\beta \leq \alpha$, the $\beta$th layer of $R_{\alpha + 1}$ is countably generated, while its ($\alpha +1$)th layer, $R_{\alpha + 1}/\psi\iota_\alpha(I)$, is isomorphic to $K$. 

\medskip
If $\alpha < \omega_1$ is a limit ordinal, then $\alpha = \sup_{n < \omega} \beta_n$ for a strictly increasing chain of countable ordinals $( \beta_n \mid n < \omega )$. We fix one such chain with $\beta_0 > 0$. By the induction hypothesis, we have the following commutative diagram, where $\nu_{\beta_n}$ ($n < \omega$) are $K$-algebra embeddings, $f_{\beta_{n+1},\beta_n}$, $g_{\beta_{n+1},\beta_n}$ ($n < \omega$) are non-unital $K$-algebra embeddings, $g_{\beta_{n+1},\beta_n}$ splits, and the complement of $\im{g_n}$ in $K^\omega$ is isomorphic to $K^\omega$ ($n < \omega$):  

\medskip
$$\begin{CD}
\dots @>>> R_{\beta_n} @>{f_{\beta_{n+1},\beta_n}}>> {R_{\beta_{n+1}}} @>>> \dots \\
@. @V{\nu_{\beta_n}}VV @V{\nu_{\beta_{n+1}}}VV @. \\
\dots @>>> {K^\omega} @>{g_{\beta_{n+1}, \beta_n}}>> {K^\omega} @>>> \dots \\
\end{CD}$$

\medskip
Let $T_\alpha = \varinjlim_{n < \omega} R_{\beta_n}$. Since the direct limit of the split embeddings in the bottom row is isomorphic to $\bigoplus_{n < \omega} K^\omega$, we obtain a non-unital $K$-algebra embedding $\varinjlim_{n < \omega} \nu_{\beta_n} : T_\alpha \to \bigoplus_{n < \omega} K^\omega$ which can be extended to a $K$-algebra homomorphism $\nu_\alpha : R_\alpha \to K^\omega$ by the same procedure as in the induction step above. 

The direct limit of the non-unital $K$-algebra embeddings in the top row yields $f_{\alpha , \beta_n} : R_{\beta_n} \hookrightarrow T_\alpha \hookrightarrow R_\alpha$ for each $n < \omega$. Moreover, we have the split non-unital $K$-algebra embeddings $g_{\alpha , \beta_n} : K^\omega \to K^\omega$ that make the corresponding squares for $\beta_n$ and $\alpha$ commute. 

Since for each $\beta < \alpha$, there exists $n < \omega$ such that $\beta < \beta_n$, we can define $f_{\alpha \beta} = f_{\alpha, \beta_n} f_{\beta_n, \beta}$ and $g_{\alpha \beta} = g_{\alpha, \beta_n} g_{\beta_n, \beta}$. Then also the squares for $\beta$ and $\alpha$ commute, by the induction premise. Finally, $R_\alpha$ is a semiartinian von Neumann regular ring with pfa of Loewy length $\alpha + 1$, all of whose layers are countably generated.     

By Lemma \ref{hereditary}, $R_\alpha$ is hereditary for each $0 < \alpha < \omega_1$. The corresponding picture for $R_\alpha$ is as follows:  

\medskip

\begin{table}[!ht]
\begin{tabular}{|l|lllllll}
\cline{1-1}
\cellcolor[HTML]{34CDF9}{$\mathbf{L_\alpha = K}$}              &  & \cellcolor[HTML]{FFFFFF}\textbf{}        & \cellcolor[HTML]{FFFFFF}\textbf{}    & \cellcolor[HTML]{34CDF9}\textbf{\,\, K}             & 
\cellcolor[HTML]{FFFFFF}             &                                      &                                      \\ \cline{1-1}
\cellcolor[HTML]{96FFFB}{...} &  & \cellcolor[HTML]{96FFFB}{...} & \cellcolor[HTML]{96FFFB}{...} & \cellcolor[HTML]{96FFFB}{...} & \cellcolor[HTML]{96FFFB}\textbf{...} & \cellcolor[HTML]{96FFFB}\textbf{...} \\ \cline{1-1} 
\cellcolor[HTML]{38FFF8}{$\mathbf{L_\beta = K^{(\omega)}}$} &  & \cellcolor[HTML]{38FFF8}{$\mathbf{K \, \oplus}$}  & \cellcolor[HTML]{38FFF8}\textbf{...} & \cellcolor[HTML]{38FFF8}{$\mathbf{\oplus \, K \, \oplus}$} & \cellcolor[HTML]{38FFF8}\textbf{...} & \cellcolor[HTML]{38FFF8}\textbf{...} \\ \cline{1-1}
\cellcolor[HTML]{68CBD0}{...} &  & \cellcolor[HTML]{68CBD0}{...} & \cellcolor[HTML]{68CBD0}{...} & \cellcolor[HTML]{68CBD0}{...} & \cellcolor[HTML]{68CBD0}\textbf{...} & \cellcolor[HTML]{68CBD0}\textbf{...} \\ \cline{1-1}
\cellcolor[HTML]{00D2CB}{$\mathbf{L_1 = K^{(\omega)}}$}  &  & \cellcolor[HTML]{00D2CB}{$\mathbf{K \, \oplus}$} & \cellcolor[HTML]{00D2CB}{...} & \cellcolor[HTML]{00D2CB}{$\mathbf{\oplus \, K \, \oplus}$ ...} & \cellcolor[HTML]{00D2CB}\textbf{...} & \cellcolor[HTML]{00D2CB}\textbf{...} \\ \cline{1-1}
\cellcolor[HTML]{329A9D}{$\mathbf{L_0 = K^{(\omega)}}$}  &  & \cellcolor[HTML]{329A9D}{$\mathbf{K \, \oplus}$}  & \cellcolor[HTML]{329A9D}{...} & \cellcolor[HTML]{329A9D}{$\mathbf{\oplus \, K \, \oplus}$ ...} & \cellcolor[HTML]{329A9D}\textbf{...} & \cellcolor[HTML]{329A9D}\textbf{...} \\ \cline{1-1}
\end{tabular}
\end{table}

\end{example}

\eject 

\medskip
Finally, we recall that in the hereditary setting of Lemma \ref{hereditary}(2), there is a single short exact sequence that tests for weak $R$-projectivity of any module $M$: 

\begin{lemma}\label{test} Let $R$ be a semiartinian von Neumann regular ring with pfa. Assume $\sigma$ is a countable ordinal and all proper layers of $R$ are countably generated. 

Then there exist a module $B$ which is a countable direct product of certain ideals of $R$, an injective module $N$, and an epimorphism $\pi : B \to N$, such that the following are equivalent for a module $M$:

\begin{enumerate}
\item $M$ is weakly $R$-projective.
\item The homomorphism $\Hom RM{\pi} : \Hom RMB \to \Hom RMN$ is surjective. 
\end{enumerate}  
\end{lemma}
\begin{proof}
This follows by \cite[4.2]{T2}.
\end{proof}

\begin{remark}\label{finer} By \cite[4.1]{T2}, the epimorphism $\pi$ is the product of restrictions of the layer epimorphisms $\pi_\alpha$ ($\alpha \leq \sigma$) to the right ideals $N_{\alpha , F}$, where $F$ runs over all finite subsets of $\lambda _\alpha$, $S_\alpha \subseteq N_{\alpha , F} \subseteq S_{\alpha + 1}$, and $N_{\alpha , F}/S_\alpha \cong \bigoplus_{\beta \in F} M_{n_{\alpha\beta}}(K_{\alpha\beta})$ is an injective module, cf.\ Theorem \ref{RTZth}. 

In particular, \cite[4.2]{T2} implies that if $M$ is not weakly $R$-projective, then there exist an $\alpha \leq \sigma$ and a finite subset $F$ of $\lambda _\alpha$ such that the homomorphism $\Hom RM{\pi_\alpha \restriction N_{\alpha , F }} : \Hom RM{N_{\alpha , F}} \to \Hom RM{N_{\alpha , F}/S_\alpha}$ is not surjective.  
\end{remark}

\section{Weak diamond and weak projectivity}

By Lemma \ref{count}(2), if $R$ is a semiartinian von Neumann regular ring with pfa, then the notions of a weak $R$-projective, $R$-projective, and projective module coincide for any countably generated module. In contrast with Lemma \ref{sup}, we will show that in the extension of ZFC + CH where the prediction principle $\Phi$ (Weak Diamond) holds, these notions coincide for \emph{arbitrary} modules, provided that condition (2) of Lemma \ref{hereditary} holds (whence $R$ is hereditary) and $\card{R} \leq \aleph_1$. 

To simplify our notation, we introduce the following definition: 
   
\begin{definition}\label{small} Let $R$ be a semiartinian von Neumann regular ring with pfa. Then $R$ is \emph{small}, if $\card{R} \leq \aleph_1$, $\sigma$ is a countable ordinal, and all proper layers of $R$ are countably generated. 
\end{definition}

Notice that our notion of smallness is more general than the one in \cite[Definition 4.3]{T2} which required the ordinal $\sigma$ to be finite rather than countable.

\medskip
Before introducing the Weak Diamond Principle, we need to recall several basic set-theoretic notions: 

Let $\kappa$ be a regular uncountable cardinal. Let $A$ be a set of cardinality $\leq \kappa$. An increasing continuous chain, $( A_\gamma \mid \gamma < \kappa )$, consisting of subsets of $A$ of cardinality $< \kappa$, is a \emph{$\kappa$-filtration of the set} $A$ in case its union is $A$. If $A$ is moreover a module, we will also use the term \emph{$\kappa$-filtration of the module} $A$, which has the additional assumptions that all the $A_\gamma$ ($\gamma < \kappa$) are submodules of $A$, and $A_0 = 0$. 

A subset $C \subseteq \kappa$ is a \emph{club} in $\kappa$, if $C$ is unbounded (i.e., $\sup C = \kappa$) and closed (i.e. for each $D \subseteq C$, if $s = \sup D < \kappa$ then $s \in C$). A subset
$E \subseteq \kappa$ is \emph{stationary} in $\kappa$, if $E \cap C \neq \emptyset$ for each club $C$ in $\kappa$. 
 
Now, we can introduce the Weak Diamond Principle, $\Phi$. We will use it in the following form presented in \cite[Lemma VI.1.7]{EM} and \cite[Theorem 18.12]{GT}: 

\medskip
($\mathbf{\Phi}$) \, Let $\kappa$ be a regular uncountable cardinal and $E$ a stationary subset in $\kappa$. Let $A$ and $B$ be sets of cardinality $\leq \kappa$. Let $( A_\gamma \mid \gamma < \kappa )$ be a $\kappa$-filtration of $A$, and $( B_\gamma \mid \gamma < \kappa )$ a $\kappa$-filtration of $B$. For each $\gamma \in E$, let $c_\gamma: {}^{A_\gamma}B_\gamma \to 2$. Then there exists a function $c : E \to 2$, such that for each $x \in {}^{A}B$, the set $E(x) = \{ \gamma \in E \mid x \restriction A_\gamma \in {}^{A_\gamma}B_\gamma \hbox{ and } c (\gamma) = c_\gamma (x \restriction A_\gamma) \}$ is stationary in $\kappa$.

\medskip
The Weak Diamond $\Phi$ is easily seen to be a consequence the better known (and stronger) Jensen's Diamond $\diamondsuit$,  which in turn is a consequence of G\"{o}del's Axiom of Constructibility, and hence is consistent with ZFC + GCH. We refer to \cite[\S VI.1]{EM} and \cite[\S 18]{GT} for more details.

\medskip
We arrive at the promised generalization of \cite[4.4]{T2} requiring only the Weak Diamond (rather than the stronger Jensen's Diamond), and allowing for arbitrary countable Loewy length of $R$:   

\begin{theorem}\label{weak} Assume CH + $\Phi$. Let $R$ be small. Then all weakly $R$-projective modules, and hence all $R$-projective modules, are projective.  
\end{theorem}
\begin{proof} Let $M$ be a weakly $R$-projective module. By induction on the minimal number of generators, $\kappa$, of $M$, we will prove that $M$ is projective. For $\kappa \leq \aleph_0$, the result follows by Lemma \ref{count}(2). If $\kappa$ is a singular cardinal, then we use the fact that the class of all weakly $R$-projective modules is closed under submodules (Lemma \ref{count}(1)) and apply Shelah's Singular Compactness Theorem (e.g., in the version from \cite[7.9]{GT}, see also \cite[XII.1.14]{EM}).   

Assume that $\kappa$ is a regular uncountable cardinal. Let $A = \{ m_\gamma \mid \gamma < \kappa \}$ be a minimal set of $R$-generators of $M$. For each $\gamma < \kappa$, let $A_\gamma = \{ m_\delta \mid \delta < \gamma \}$. Let $M_\gamma$ be the submodule of $M$ generated by $A_\gamma$. By the inductive premise, $M_\gamma$ is projective, and $\mathcal M = (M_\gamma \mid \gamma < \kappa )$ is a $\kappa$-filtration of the module $M$. Possibly skipping some terms of $\mathcal M$, we can w.l.o.g.\ assume that $\mathcal M$ has the following property for each $\gamma < \kappa$: if $M_{\delta}/M_\gamma$ is not weakly $R$-projective for some $\gamma < \delta < \kappa$, then also $M_{\gamma + 1}/M_\gamma$ is not weakly $R$-projective. 

\medskip
($\dagger$) Let $E$ be the set of all $\gamma < \kappa$ such that $M_{\gamma + 1}/M_\gamma$ is not weakly $R$-projective. Also, let $\pi : B \to N$ be the epimorphism from Lemma \ref{test}. By that Lemma, for each $\gamma \in E$, we can choose an $h_\gamma \in \Hom R{M_{\gamma + 1}/M_\gamma}N$ such that $h_\gamma$ does not factorize through $\pi$. 

\medskip 
We claim that $E$ is not stationary in $\kappa$. Assume this claim is not true. Note that CH implies $\card{B} \leq \aleph_1 \leq \kappa$, so we can fix a $\kappa$-filtration of the set $B$, $( B_\gamma \mid \gamma < \kappa )$. 

Let $\gamma < \kappa$. For each $g \in \Hom R{M_\gamma}N$, we choose $g^+ \in \Hom R{M_{\gamma +1}}N$ such that $g^+ \restriction M_\gamma = g$. This is possible because $N$ is an injective module.
Further, for each $f \in \Hom R{M_\gamma}B$, we choose $f^e \in \Hom R{M_{\gamma +1}}B$ such that $\pi f^e = (\pi f)^+$. This is possible since $M_{\gamma + 1}$ is projective. Notice that 
$\pi(f^e \restriction M_\gamma) = (\pi f)^+ \restriction M_\gamma = \pi f$, so $\delta_f := f^e \restriction M_\gamma - f \in \Hom R{M_\gamma}K$ where $K = \Ker \pi$.   

\medskip
For each $\gamma \in E$, we define $c_\gamma : {}^{A_\gamma}B_\gamma \to 2$ as follows: If $x : A_\gamma \to B_\gamma$ is a restriction of a (necessarily unique) morphism $f \in \Hom R{M_\gamma}B$ such that the morphism $\delta_f = f^e \restriction M_\gamma - f$ can be extended to a morphism from $\Hom R{M_{\gamma + 1}}K$, then we put $c_\gamma (x) = 1$. Otherwise, we let $c_\gamma (x) = 0$.      

In this setting, $\Phi$ yields a function $c : E \to 2$ such that for each $x \in {}^{A}B$, the set $E(x) = \{ \gamma \in E \mid x \restriction A_\gamma \in {}^{A_\gamma}B_\gamma \hbox{ and } c (\gamma) = c_\gamma (x \restriction A_\gamma) \}$ is stationary in $\kappa$. We will use $c$ to define a morphism $g \in \Hom RMN$ as follows:

\medskip
By induction on $\gamma < \kappa$, we define a sequence $( g_\gamma \mid \gamma < \kappa )$ such that $g_\gamma \in \Hom R{M_\gamma}N$. First, $g_0 = 0$. If $\gamma < \kappa$, and $g_\gamma$ is already defined, we distinguish two cases: 
  
(I) $\gamma \notin E$ or $c(\gamma) = 0$. In this case, we put $g_{\gamma + 1} = (g_\gamma)^+$.

(II) $\gamma \in E$ and $c(\gamma) = 1$. In this case, we let $g_{\gamma + 1} = (g_\gamma)^+ + h_\gamma \rho_\gamma$, where $\rho_\gamma : M_{\gamma + 1} \to M_{\gamma + 1}/M_\gamma$ is the canonical projection modulo $M_\gamma$.

Notice that in both cases $g_{\gamma + 1} \restriction M_\gamma = g_\gamma$. We let $g_\gamma = \bigcup_{\delta < \gamma} g_\delta$ in case $\gamma < \kappa$ is a limit ordinal. Then $g = \bigcup_{\gamma < \kappa} g_\gamma \in \Hom RMN$.

\medskip
Since $M$ is weakly $R$-projective, there exists $f \in \Hom RMB$ such that $g = \pi f$. By $\Phi$, the set $E(f \restriction A) = \{ \gamma \in E \mid f \restriction A_\gamma \in {}^{A_\gamma}B_\gamma \hbox{ and } c (\gamma) = c_\gamma (f \restriction A_\gamma) \}$ is stationary in $\kappa$. Let $\gamma \in E(f \restriction A)$.    

Assume that $c(\gamma) = 0$. Then we are in case (I), so $\pi (f \restriction M_{\gamma + 1}) = g_{\gamma + 1} = (g_\gamma)^+ = ( \pi f \restriction M_\gamma )^+ = \pi (f \restriction M_\gamma )^e$. Then the morphism $(f \restriction M_\gamma )^e - f \restriction M_{\gamma + 1} \in \Hom R{M_{\gamma + 1}}K$ is an extension of $\delta _{f \restriction M_\gamma} = (f \restriction M_\gamma) ^e \restriction M_\gamma - f \restriction M_\gamma$ to $M_{\gamma + 1}$, in contradiction with $c_\gamma (f \restriction A_\gamma) = c(\gamma) = 0$. 

So necessarily $c(\gamma) = 1$, and we are in case (II), whence $g_{\gamma + 1} = (g_\gamma)^+ + h_\gamma \rho_\gamma$. As $c_\gamma (f \restriction A_\gamma) = c(\gamma) = 1$, the morphism $\delta_{f \restriction M_\gamma} =  (f \restriction M_\gamma ) ^e \restriction M_\gamma -  f \restriction M_\gamma$ can be extended to a morphism $\Delta _{f \restriction M_\gamma} \in \Hom R{M_{\gamma + 1}}K$. 

Again, $g_{\gamma + 1} = \pi (f \restriction M_{\gamma + 1})$ and $(g_\gamma)^+ =  \pi (f \restriction M_\gamma)^e$, so 
$$h_\gamma \rho_\gamma = \pi (f \restriction M_{\gamma + 1} - (f \restriction M_\gamma)^e) = \pi u_\gamma ,$$
where $u_\gamma = f \restriction M_{\gamma + 1} - (f \restriction M_\gamma)^e + \Delta _{f \restriction M_\gamma} \in \Hom R{M_{\gamma + 1}}B$. However, $u_\gamma \restriction M_\gamma = f \restriction M_{\gamma} - (f \restriction M_\gamma)^e \restriction M_\gamma + \delta_{f \restriction M_\gamma} = 0$. Hence $u_\gamma$ factorizes through $\rho_\gamma$ by some $v_\gamma \in \Hom R{M_{\gamma + 1}/M_\gamma}B$, that is, $u_\gamma = v_\gamma \rho_\gamma$.

Thus $h_\gamma \rho_\gamma = \pi v_\gamma \rho_\gamma$. Since $\rho_\gamma$ is surjective, $h_\gamma = \pi v_\gamma$, in contradiction with our choice of $h_\gamma$.             

\medskip  
This proves our claim about the set $E$. So there is a club $C$ in $\kappa$ such that $C \cap E = \emptyset$. Let $z : \kappa \to \kappa$ be a strictly increasing continuous function whose image is $C$. For each $\gamma < \kappa$, let $N_\gamma = M_{z(\gamma )}$. Then $(N_\gamma \mid \gamma < \kappa )$ is a $\kappa$-filtration of the module $M$ such that $N_{\gamma + 1}/N_\gamma$ is weakly $R$-projective for each $\gamma < \kappa$. By the inductive premise, $N_{\gamma + 1}/N_{\gamma}$ is projective, hence $N_{\gamma + 1} = N_\gamma \oplus P_\gamma$ for a projective module $P_\gamma$, for each $\gamma < \kappa$. We conclude that $M = N_0 \oplus \bigoplus_{\gamma < \kappa} P_\gamma$ is projective. This finishes the inductive step for the case when $\kappa$ is a regular uncountable cardinal.  
\end{proof}           

\begin{remark}\label{JSar} As observed by Jan \v Saroch, Remark \ref{finer} makes it possible to prove Theorem \ref{weak} without the assumption of CH, i.e., assuming only $\Phi$. (That is indeed a weaker assumption, since unlike $\diamondsuit$, the $\Phi$ does not imply CH.) 

The only modification needed concerns the set $E$ of all $\gamma < \kappa$ such that $M_{\gamma + 1}/M_\gamma$ is not weakly $R$-projective, and the morphisms $h_\gamma$ ($\gamma \in E$) defined in part ($\dagger$) of the proof of Theorem \ref{weak}: 

For each $\gamma \in E$, Remark \ref{finer} yields an $\alpha _\gamma \leq \sigma$, a finite subset $F_\gamma$ of $\lambda _{\alpha_\gamma}$, and an $h_\gamma \in \Hom R{M_{\gamma + 1}/M_\gamma}{N_{\alpha _\gamma , F _\gamma }/S_\alpha}$ which does not factorize through $\pi_\alpha \restriction N_{\alpha _\gamma , F _\gamma }$. For each $\alpha \leq \sigma$ and each finite subset $F$ of $\lambda _{\alpha}$, let $E_{\alpha , F } = \{ \gamma \in E \mid \alpha _\gamma = \alpha \, \& \, F_\gamma = F \}$. Then $E = \bigcup_{\alpha , F } E_{\alpha , F }$. Notice that since $R$ is small, the set of all such pairs $(\alpha , F )$ is countable. 

Thus, if $E$ is stationary, then so is one of the $E_{\alpha , F }$ (see e.g.\ \cite[II.4.3]{EM}), say $E_{\alpha ^\prime, F ^\prime}$. The proof of a contradiction with the stationarity of $E$ then proceeds as that of Theorem \ref{weak} in the parts following $(\dagger$), but for $B = N_{\alpha ^\prime , F ^\prime}$, $N = N_{\alpha ^\prime, F ^\prime}/S_{\alpha ^\prime}$, $\pi = \pi_{\alpha ^\prime} \restriction N_{\alpha ^\prime , F ^\prime }$, and $E = E_{\alpha ^\prime, F ^\prime}$. The point is that in this setting, the cardinality of $B$ is $\leq \card{R} \leq \aleph_1$ even without the assumption of CH. 
\end{remark}

Combining Lemma \ref{sup} and Theorem \ref{weak}, we obtain

\begin{corollary}\label{indep} Let $R$ be small. Then the assertion {\lq}All weakly $R$-projective modules are projective{\rq} is independent of ZFC + GCH.
\end{corollary}

\begin{acknowledgment}
The author thanks Kate\v{r}ina Fukov\'{a} for valuable comments on an earlier draft of this paper.   
\end{acknowledgment}

\eject

\end{document}